\documentclass[numreferences]{kluwer}
\usepackage{amssymb,graphicx}

\newtheorem{theo}{Theorem}
\newtheorem{defn}{Definition}

\begin{document}
\begin{article}
\begin{opening}
\markboth{G. A. Panopoulos, Z. A. Anastassi and T. E. Simos}
{Two Optimized Symmetric Eight-Step Implicit Methods for IVPs with Oscillating Solutions}

\title
{{Two Optimized Symmetric Eight-Step Implicit Methods for Initial-Value Problems with Oscillating Solutions}}

\author{G.A. \surname{Panopoulos}\thanks{e-mail: gpanop@uop.gr}}
\author{Z.A. \surname{Anastassi}\thanks{e-mail: zackanas@uop.gr}}
\author{T.E. \surname{Simos}\thanks{Highly Cited Researcher, Active Member of the European Academy of Sciences and Arts, Address: Dr. T.E. Simos, 26 Menelaou Street, Amfithea - Paleon Faliron, GR-175 64 Athens, GREECE, Tel: 0030 210 94 20 091, e-mail: tsimos.conf@gmail.com, tsimos@mail.ariadne-t.gr}}
 \institute{Laboratory of Computer Sciences, Department of Computer Science and Technology, Faculty of Sciences and Technology, University of
Peloponnese,\\ GR-22 100 Tripolis, GREECE}

\runningtitle{Optimized Symmetric 8-Step Implicit Methods for IVPs with Oscillating Solutions}
\runningauthor{G. A. Panopoulos, Z. A. Anastassi, T. E. Simos}

\begin{abstract}
In this paper we present two optimized eight-step symmetric implicit methods with phase-lag order ten and infinite (phase-fitted). The
methods are constructed to solve numerically the radial time-independent Schr\"odinger equation with the use of the Woods-Saxon potential. They can also be used to integrate related IVPs with oscillating solutions such as orbital problems. We compare the two new methods to some recently constructed optimized methods from the literature. We measure the efficiency of the methods and conclude that the new method with infinite order of phase-lag is the most efficient of all the compared methods and for all the problems solved.
\end{abstract}


\classification{PACS}{0.260, 95.10.E}
\end{opening}

\section{Introduction}
\label{Intro}

The radial Schr\"{o}dinger equation can be written as:

\begin{equation}
\label{Schrodinger}
    y''(r) = \left( \frac{l(l+1)}{r^{2}}+V(r)-E \right) y(r)
\end{equation}

\noindent where $\frac{l(l+1)}{r^{2}}$ is the \textit{centrifugal potential}, $V(r)$ is the \textit{potential}, $E$ is
the \textit{energy} and $W(r) = \frac{l(l+1)}{r^{2}} + V(r)$ is the \textit{effective potential}. It is valid that
${\mathop {\lim} \limits_{r \to \infty}} V(r) = 0$ and therefore ${\mathop {\lim} \limits_{r \to \infty}} W(r) = 0$.

We consider $E>0$ and divide $[0,\infty)$ into subintervals $[a_{i},b_{i}]$ so that $W(r)$ is a constant with value
${\mathop{W_{i}}\limits^{\_}}$. After this the problem (\ref{Schrodinger}) can be expressed by the approximation:

\begin{equation}
\begin{array}{l}
\label{Schrodinger_simpl}
y''_{i} = ({\mathop{W}\limits^{\_}} - E)\,y_{i},\\
\end{array}
\end{equation}
whose solution is:
\begin{equation}
\begin{array}{l}
\label{Schrodinger_simplified}
y_{i}(r) = A_{i}\,\exp{\left(\sqrt{{\mathop{W}\limits^{\_}}-E}\,r\right)} +
B_{i}\,\exp{\left(-\sqrt{{\mathop{W}\limits^{\_}}-E}\,r\right)}, \\
 A_{i},\,B_{i}\,\in {\mathbb R}.
\end{array}
\end{equation}

Many numerical methods have been developed for the efficient solution of the Schr\"odinger equation and related
problems. For example Raptis and Allison have developed a two-step exponentially-fitted method of order four in
\cite{rallison}. More recently Kalogiratou and Simos have constructed a two-step P-stable exponentially-fitted method of
order four in \cite{kalogiratou}.

Some other notable multistep methods for the numerical solution of oscillating IVPs have been developed by Chawla and
Rao in \cite{chawla}, who produced a three-stage, two-Step P-stable method with minimal phase-lag and order six and by
Henrici in \cite{henrici}, who produced a four-step symmetric method of order six. Also Anastassi and Simos have developed trigonometrically fitted six-step symmetric methods in \cite{anastassi_simos_p8} and a six-step P-stable trigonometrically-fitted method in \cite{anastassi_simos_p11} and Panopoulos, Anastassi and Simos have constructed two optimized eight-step symmetric methods.

Also some research work in numerical methods can be found in \cite{royal}-\cite{jnaiam3_11}.

\section{Phase-lag analysis of symmetric multistep methods}

For the numerical solution of the initial value problem

\begin{equation}
\label{ivp_definition}
    y'' = f(x,y)
\end{equation}

\noindent multistep methods of the form

\begin{equation}
\label{multistep_definition}
    \sum\limits_{i=0}^{m}{a_{i}y_{n+i}} = h^{2}\sum\limits_{i=0}^{m}{b_{i}f(x_{n+i},y_{n+i})}
\end{equation}

with $m$ steps can be used over the equally spaced intervals $\left\{x_{i}\right\}^{m}_{i=0} \in [a,b]$ and
$h=|x_{i+1}-x_{i}|$, \, $i=0(1)m-1$.

If the method is symmetric then $a_i=a_{m-i}$ and $b_i=b_{m-i}$, \, $i=0(1)\lfloor \frac{m}{2} \rfloor$.

Method (\ref{multistep_definition}) is associated with the operator

\begin{eqnarray}
\label{exp_operator} L(x) = \sum\limits_{i=0}^{m}{a_{i}u(x+ih)} - h^{2}\sum\limits_{i=0}^{m}{b_{i}u''(x+ih)}
\end{eqnarray}

\noindent where $u \in C^2$.

\begin{defn}
\label{defn_exp1} \emph{The multistep method (\ref{multistep_definition}) is called algebraic of order $p$ if the associated
linear operator $L$ vanishes for any linear combination of the linearly independent functions $\,1,\, x,\,x^2,\, \ldots
,\, x^{p+1}$.}
\end{defn}

When a symmetric $2k$-step method, that is for $i=-k(1)k$, is applied to the scalar test equation

\begin{equation}
\label{stab_eq} y''=-\omega^2 y
\end{equation}

a difference equation of the form
\begin{eqnarray}
\label{phl_multi_de}
\nonumber A_{k} (v)y_{n + k} + ... + A_{1} (v)y_{n + 1} + A_{0} (v)y_{n}\\
+ A_{1}(v)y_{n - 1} + ... + A_{k} (v)y_{n - k} = 0
\end{eqnarray}

\noindent is obtained, where $v = \omega h$, $h$ is the step length and $A_{0} (v)$, $A_{1} (v),\ldots$, $ A_{k} (v)$
are polynomials of $v$.

The characteristic equation associated with (\ref{phl_multi_de}) is
\begin{eqnarray}
\label{phl_multi_ce}
A_{k} (v)s^{k} + ... + A_{1} (v)s + A_{0} (v) + A_{1} (v)s^{ - 1} + ... + A_{k} (v)s^{ - k} = 0
\end{eqnarray}
\\
From Lambert and Watson (1976) we have the following definitions:

\begin{defn}
\label{defn_exp2} \emph{A symmetric $2k$-step method with characteristic equation given by (\ref{phl_multi_ce}) is said to have an interval of periodicity $(0,v_0^2)$ if, for all $v\in(0,v_0^2),$ the roots $s_i,i=1(1)2k$ of Eq. (\ref{phl_multi_ce}) satisfy:
}
{\begin{equation}
s_1=e^{i\theta(v)},\, s_2=e^{-i\theta(v)},\,\, and \,\, |s_i| \leq 1, \, i=3(1)2k
\end{equation}}{where $\theta(v)$ is a real function of v.}
\end{defn}

\begin{defn}
\label{defn_exp3} \emph{For any method corresponding to the characteristic equation (\ref{phl_multi_ce}) the phase-lag is defined as the leading term in the expansion of }

{\begin{equation}
 t=v-\theta(v)
\end{equation}}{Then if the quantity $t=O(v^{q+1})\,\,as\,\,{v \to \infty}, $ the order phase-lag is q.}
\end{defn}

\begin{theo}
\emph{\cite{royal}} The symmetric $2k$-step method with characteristic equation given by (\ref{phl_multi_ce}) has
phase-lag order $q$ and phase-lag constant $c$ given by

\begin{equation}
\label{phl_multi_defn} - c v ^{q + 2} + O(v^{q + 4}) = {\frac{{2A_{k} (v)\cos (k v ) + ... + 2A_{j} (v)\cos (j v ) + ...
+ A_{0} (v)}}{{2k^{2}A_{k} (v) + ... + 2j^{2}A_{j} (v) + ... + 2A_{1} (v)}}}
\end{equation}
\end{theo}

The formula proposed from the above theorem gives us a direct method to calculate the phase-lag of any symmetric $2k$-
step method.

\ In our case, the symmetric 8-step method has
phase-lag order $q$ and phase-lag constant $c$ given by:
\begin{equation}
\begin{array}{l}{- c v ^{q + 2} + O(v^{q + 4}) =}\\
{\frac{{2A_{4} (v)\cos (4 v ) + 2A_{3} (v)\cos (3 v )  + 2A_{2} (v)\cos (2 v ) + 2A_{1} (v)
\cos ( v )+ A_{0} (v)}}{{32A_{4} (v) + 18A_{3} (v) + 8A_{2} (v) + 2A_{1} (v)}}}
\end{array}
\end{equation}
\section{Construction of the new optimized multistep methods}
\label{Construction}

We consider the eight-step symmetric implicit methods of the form:

\begin{equation}
\begin{array}{c} \label{table_qt88}
y_{{4}} = -y_{{-4}} -a_{{3}}(y_{{3}}+y_{{-3}}) -a_{{2}}(y_{{2}}+y_{{-2}}) -a_{{1}}(y_{{1}}+y_{{-1}})\\
+{h}^{2}\left(b_{{4}}(f_{{4}}+f_{{-4}})+b_{{3}}(f_{{3}}+f_{{-3}}) +b_{{2}}(f_{{2}}+f_{{-2}}) +b_{{1}}(f_{{1}}+f_{{-1}}) +b_{{0}}f_{{0}}\right)
\end{array}
\end{equation}

\noindent where $a_{3}=-2, \quad a_{2}=2, \quad a_{1}=-1$,\\
 $y_{i} = y(x+ih)$ and $f_{i} = f(x+ih,y(x+ih))$

\subsection{First optimized method with infinite order of \\ phase-lag (phase-fitted)}
\label{Constr1}

We want the first method to have infinite order of phase-lag, that is the phase-lag will be nullified using $b_{4}$
coefficient.

We satisfy as many algebraic equations as possible, but we keep $b_{4}$ free. After achieving 10th algebraic order,
the coefficients now depend on
$b_{4}$:

\begin{equation}
\begin{array}{c}
b_{{0}}=70\,b_{{4}}-{\frac {12629}{3024}}, \qquad b_{{1}}=-56\,b_{{4}}+{ \frac {20483}{4032}},\\
\\
 \qquad b_{{2}}=28\,b_{{4}}-{\frac
{3937}{2016}},\qquad b_{{3}}=-8\,b_{{4}}+{\frac {17671}{12096}}\\
\end{array}
\end{equation}

and the phase-lag becomes:
\begin{equation}
\begin{array}{l}
PL = \frac{1}{1260}\,\frac{A}{B}, \qquad \mbox{where}\\
A = 24192\, \left( \cos \left( v \right)  \right) ^{4}+24192\, \left( \cos \left( v \right)  \right) ^{4}{v}^{2}b_{{4}}+17671\,
 \left( \cos \left( v \right)  \right) ^{3}{v}^{2}\\-96768\, \left( \cos \left( v
 \right)  \right) ^{3}{v}^{2}b_{{4}}
 -24192\,\left( \cos \left( v \right)  \right) ^{3}+14152\, \left( \cos \left( v
 \right)  \right) ^{2}{v}^{2}b_{{4}}\\-12096\, \left( \cos \left( v
 \right)  \right) ^{2}-11811\, \left( \cos \left( v \right)  \right) ^{2
}{v}^{2}+2109\,\cos \left( v \right){v}^{2}+\\
15120\,\cos \left( v \right) -96768\,\cos \left( v \right) {v}^{2}b_{{4}}-409\,{v}^{2}+24192\,{v}^{2}b_
{{4}}-3024 \quad \mbox{and}\\
B = {12+25\,{v}^{2}}
\end{array}
\end{equation}

so by satisfying $PL=0$, we derive
\begin{equation}
\label{new1}
\begin{array}{l}
\nonumber b_{4} = {-\frac {1}{24192}}\,{\frac {C}{D}}, \qquad \mbox{where}\\
\nonumber C=24192\, \left( \cos \left( v \right)  \right) ^{4}+\left( 17671\,{v}^{2}-24192\right)\left( \cos \left( v \right) \right) ^{3}\\-\left(12096+11811{v}^{2}\right) \left( \cos \left( v \right) \right) ^{2}+\\
\nonumber  \left(15120+2109\,{v}^{2} \right)\cos \left( v \right)
-409\,{v}^{2}-3024\\
D={{v}^{ 2} \left( \cos \left( v \right)^{4} -4 \cos \left( v \right)^{3}+6 \cos \left( v \right)^{2}-4\cos \left( v \right)+1\right) }
\end{array}
\end{equation}

\noindent where $v=\omega\,h$, $\omega$ is the frequency and $h$ is the step length used.

\subsection{Second optimized method with tenth order of \\ phase-lag}
\label{Constr2}

For this method we use all $b_{i}$ coefficients for achieving maximum algebraic order or maximum phase-lag order. After
achieving maximum algebraic order, that is ten, the coefficients become:
\begin{eqnarray}
\label{new2} b_{0}=\frac{17273}{72576}, \, b_{1}=\frac{280997}{181440}, \, b_{2}=-\frac{33961}{181440},\,
b_{3}=\frac{173531}{181440}, \, b_{4}=\frac{45767}{725760}
\end{eqnarray}\\
If we repeat the procedure of the previous section and expand phase-lag using the Taylor series, we can nullify the
leading term (that is the coefficient of $h^{10}$). However we obtain the same method as (\ref{new2}). The same method will
be produced if we attempt any combination of algebraic order and phase-lag order. This happens due to the symmetry of
the specific $a_i$.

\section{Numerical results}
\label{Numerical_results}

\subsection{The problems}

The efficiency of the two newly constructed methods will be measured through the integration of five initial value
problems with oscillating solution.

\subsubsection{Orbital Problem by Franco and Palacios}
The "almost" periodic orbital problem studied by \cite{franco} can be described by

\begin{equation}\label{ivp_franco_original}
y''+y=\epsilon\,e^{i\,\psi\,x}, \quad y(0)=1, \quad y'(0)=i, \quad y\in \mathcal{C},
\end{equation}

or equivalently by

\begin{equation}\label{ivp_franco_equivalent}
\begin{array}{l}
u''+u=\epsilon\,\cos(\psi\,x), \quad u(0)=1, \quad u'(0)=0,\\
v''+v=\epsilon\,\sin(\psi\,x), \quad v(0)=0, \quad v'(0)=1,
\end{array}
\end{equation}

where $\epsilon=0.001$ and $\psi=0.01$.

The theoretical solution of the problem (\ref{ivp_franco_original}) is given below:

\begin{equation}\label{ivp_franco_theoretical}
\begin{array}{l}
y(x)=u(x)+i\,v(x), \quad u,v\in\mathcal{R}\nonumber\\
u(x)={\frac {1-\epsilon-{\psi}^{2}}{1-{\psi}^{2}}}\,\cos(x) + {\frac
{\epsilon}{1-{\psi}^{2}}}\,\cos(\psi\,x)\\
v(x)={\frac {1-\epsilon\psi-{\psi}^{2}}{1-{\psi}^{2}}}\,\sin(x) + { \frac {\epsilon}{1-{\psi}^{2}}}\,\sin(\psi\,x)
\end{array}
\end{equation}

The system of equations (\ref{ivp_franco_equivalent}) has been solved for $x\in[0,1000\,\pi].$ The estimated frequency
is $w=1$.

\subsubsection{Inhomogeneous Equation}
\hspace{12pt} $y''=-100\,y+99\,\sin(t),\;$ with $\; y(0)=1, y'(0)=11,\; t\in[0,1000\,\pi]$.\vspace{5pt}\\
Theoretical solution: $y(t)=\sin(t)+\sin(10\,t)+\cos(10\,t)$.\\
Estimated frequency: $w=10$.

\subsubsection{Two-Body Problem}
\hspace{12pt} $y''=-\frac{y}{(y^2+z^2)^\frac{3}{2}},\; z''=-\frac{z}{(y^2+z^2)^\frac{3}{2}},\;$ with $y(0)=1,\;
y'(0)=0,\; z(0)=0,\; z'(0)=1, t\in[0,1000\,\pi].\;$ Theoretical solution: $y(t)=\cos(t)$ and $z(t)=\sin(t)$. We used the
estimation $w=\frac{1}{(y^2+z^2)^\frac{3}{4}}$ as frequency of the problem.

\subsubsection{Duffing Equation}
\hspace{12pt} $y''=-y-y^3+0.002\cos(1.01\,t),\;$ with $\; y(0)=0.200426728067, y'(0)=0,\; t\in[0,1000\,\pi]$.\vspace{5pt}\\
Theoretical solution: $y(t)=0.200179477536\,\cos(1.01\,t)\,+\,2.46946143\,\cdot\,{10}^{-4}\,\cos(3.03\,t)\,+\,3.04014\,\cdot\,{10}^{-7}\,
\cos(5.05\,t)\,+\,3.74\,\cdot\,{10}^{-10}\,\cos(7.07\,t)\,+\,...\,\,$.\\
Estimated frequency: $w=1$.

\subsubsection{The inverse resonance problem}
We will integrate problem (\ref{Schrodinger}) (where $r=x$)\,with $l=0$ at the interval $[0,15]$ using the well known Woods-Saxon
potential
\begin{eqnarray}
\label{Woods_Saxon} V(x) = \frac{u_{0}}{1+q} + \frac{u_{1}\,q}{(1+q)^2}, \quad\quad q =
\exp{\left(\frac{x-x_{0}}{a}\right)}, \quad
\mbox{where}\\
\nonumber u_{0}=-50, \quad a=0.6, \quad x_{0}=7 \quad \mbox{and} \quad u_{1}=-\frac{u_{0}}{a}
\end{eqnarray}

\noindent and with boundary condition $y(0)=0$.

\noindent The potential $V(x)$ decays more quickly than $\frac{l\,(l+1)}{x^2}$, so for large $x$ (asymptotic region) the
Schr\"{o}dinger equation (\ref{Schrodinger}) becomes

\begin{equation}
\label{Schrodinger_reduced}
    y''(x) = \left( \frac{l(l+1)}{x^{2}}-E \right) y(x)
\end{equation}

\noindent The last equation has two linearly independent solutions $k\,x\,j_{l}(k\,x)$ and\\
$k\,x\,n_{l}(k\,x)$, where $j_{l}$ and $n_{l}$ are the \textit{spherical Bessel} and \textit{Neumann} functions. When $x
\rightarrow \infty$ the solution takes the asymptotic form

\begin{equation}
\label{asymptotic_solution}
\begin{array}{l}y(x) \approx A\,k\,x\,j_{l}(k\,x) - B\,k\,x\,n_{l}(k\,x) \\
\approx D[sin(k\,x - \pi\,l/2) + \tan(\delta_{l})\,\cos{(k\,x - \pi\,l/2)}],
\end{array}
\end{equation}

\noindent where $\delta_{l}$ is called \textit{scattering phase shift} and it is given by the following expression:
\begin{equation}
\tan{(\delta_{l})} = \frac{y(x_{i})\,S(x_{i+1}) - y(x_{i+1})\,S(x_{i})} {y(x_{i+1})\,C(x_{i}) - y(x_{i})\,C(x_{i+1})},
\end{equation}

\noindent where $S(x)=k\,x\,j_{l}(k\,x)$, $C(x)=k\,x\,n_{l}(k\,x)$ and $x_{i}<x_{i+1}$ and both belong to the asymptotic
region. Given the energy we approximate the phase shift, the accurate value of which is $\pi/2$ for the above problem.

We will use three different values for the energy: i) $989.701916$ and ii) $341.495874$ and iii) $163.215341$. As for the frequency $\omega$ we
will use the suggestion of Ixaru and Rizea \cite{ix_ri}:

\begin{equation}
\omega = \cases{ \sqrt{E-50} & $x\in[0,\,6.5]$ \cr
            \sqrt{E}    &$x\in[6.5,\,15]$ \cr}
\end{equation}

\subsection{The methods}

We have used several multistep methods for the integration of the Schr\"odinger equation. These are:
\begin{itemize}
\item The new method with infinite order of phase-lag shown in (\ref{new1})
\item The new method with eighth order of phase-lag shown in (\ref{new2})
\item The P-stable method of Henrici with minimal phase-lag and order six \cite{henrici}
\item The three-stage method of Chawla and Rao of order six \cite{chawla}
\item The Classical method of Numerov
\item The P-stable exponentially-fitted method of Kalogiratou and Simos of order four \cite{kalogiratou}
\item The three-step method of Adams-Moulton
\end{itemize}

\subsection{Comparison}
\begin{figure}[htbp]
    \includegraphics[width=\textwidth]{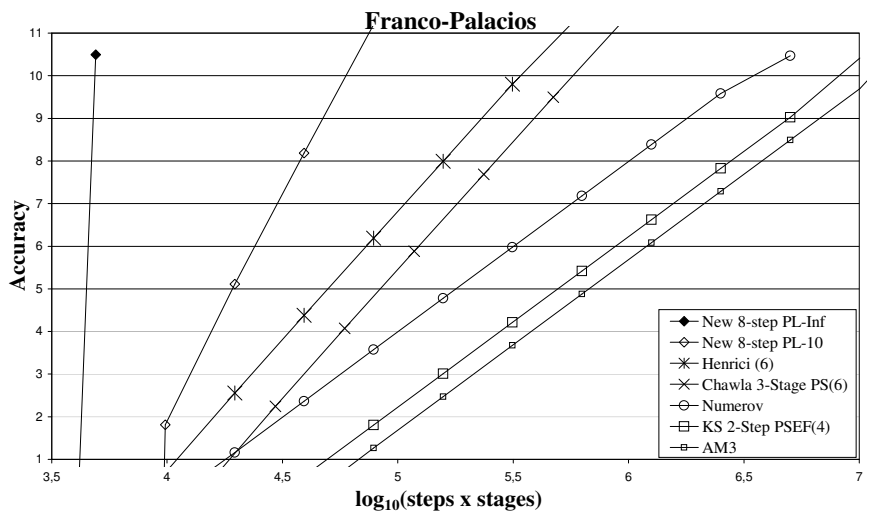}
    \caption{Efficiency for the Orbital Problem by Franco and Palacios}
    \label{fig_Orbital Problem by Franco and Palacios}
\end{figure}

\begin{figure}[htbp]
    \includegraphics[width=\textwidth]{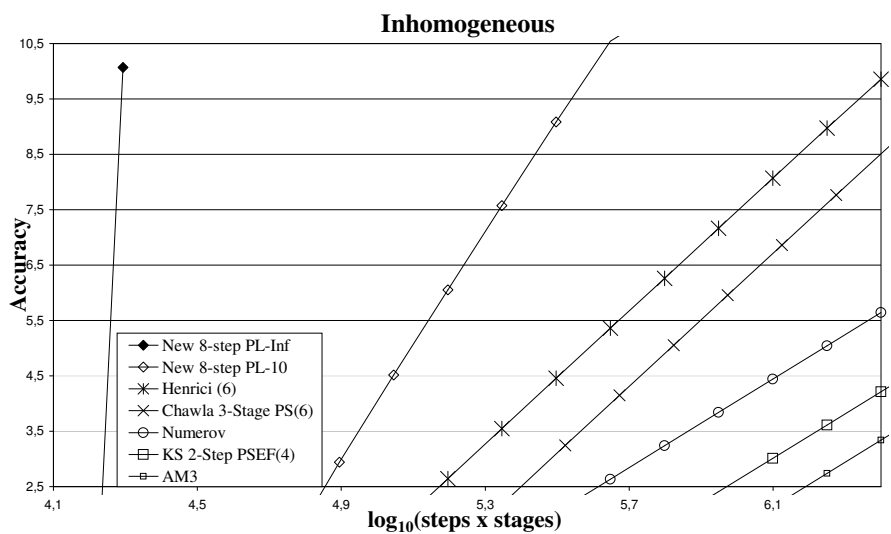}
    \caption{Efficiency for the Inhomogeneous equation}
    \label{fig_inhom}
\end{figure}

\begin{figure}[htbp]
    \includegraphics[width=\textwidth]{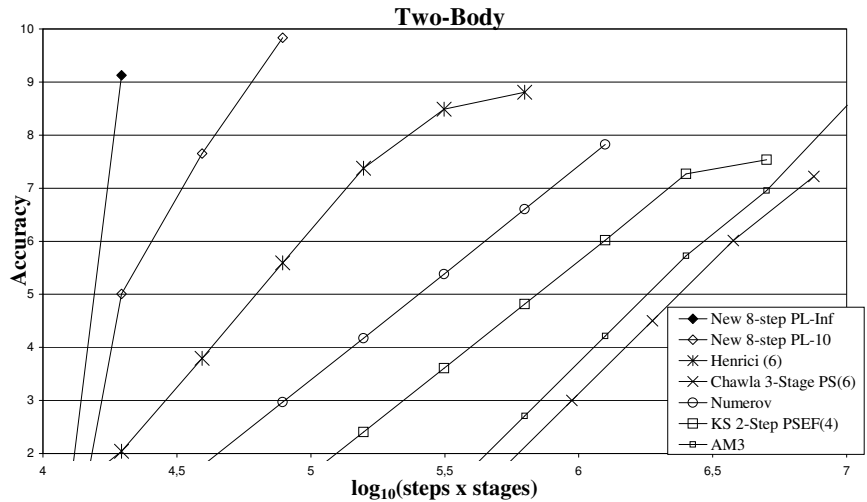}
    \caption{Efficiency for the Two-body problem}
    \label{fig_twobody}
\end{figure}

\begin{figure}[htbp]
    \includegraphics[width=\textwidth]{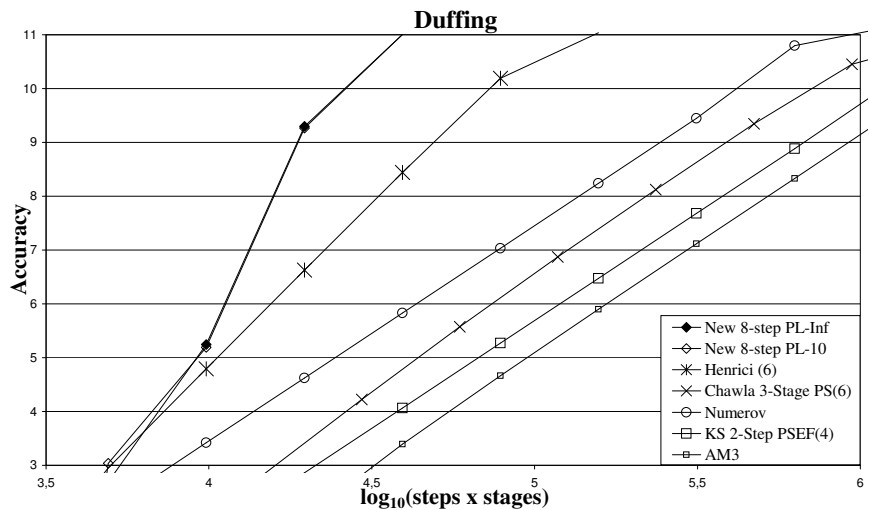}
    \caption{Efficiency for the Duffing equation}
    \label{fig_Duffing}
\end{figure}

\begin{figure}[htbp]
    \includegraphics[width=\textwidth]{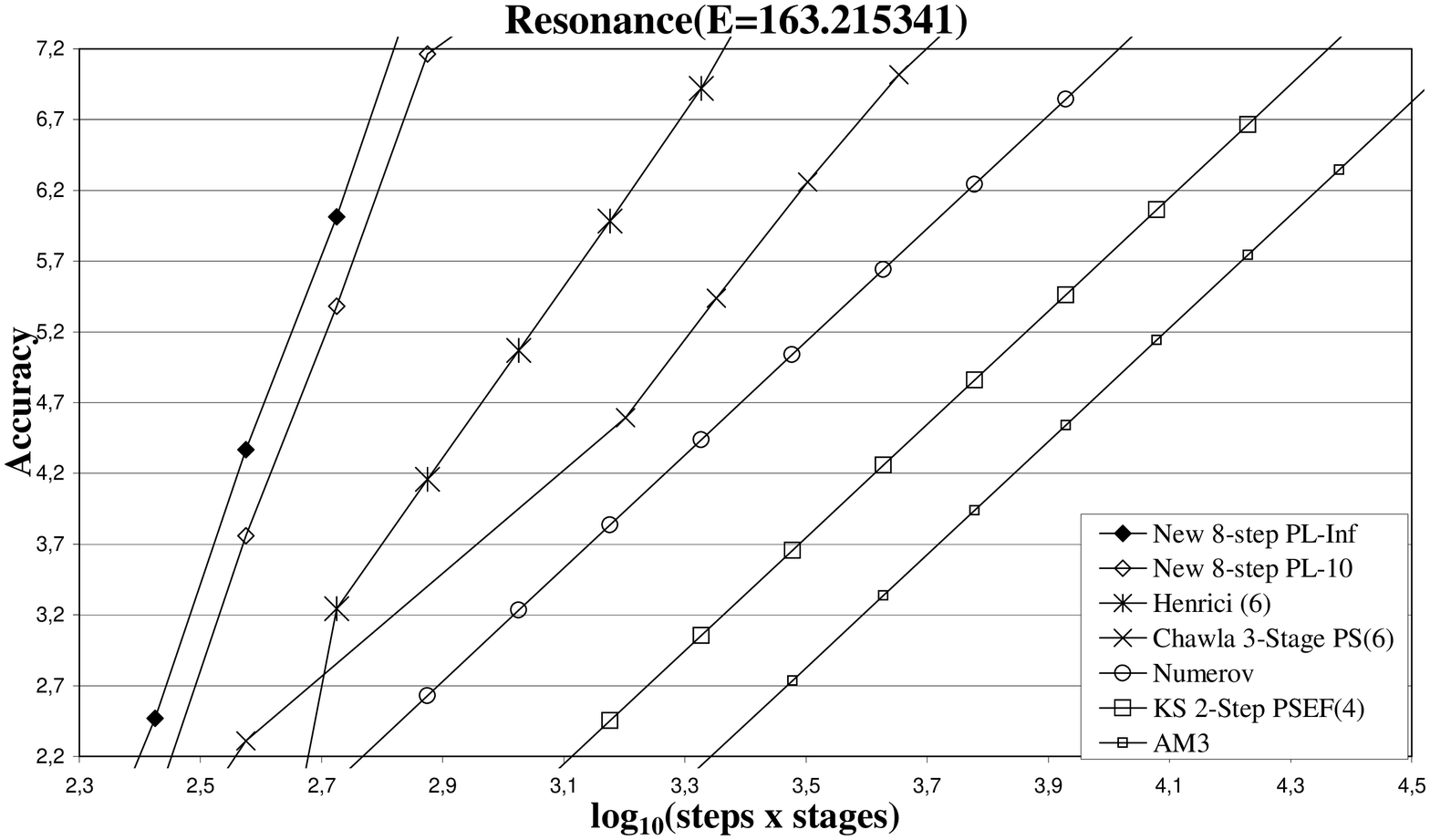}
    \caption{Efficiency for the Schr\"odinger equation using E = 163.215341}
    \label{fig_resonance_163}
\end{figure}

\begin{figure}[htbp]
    \includegraphics[width=\textwidth]{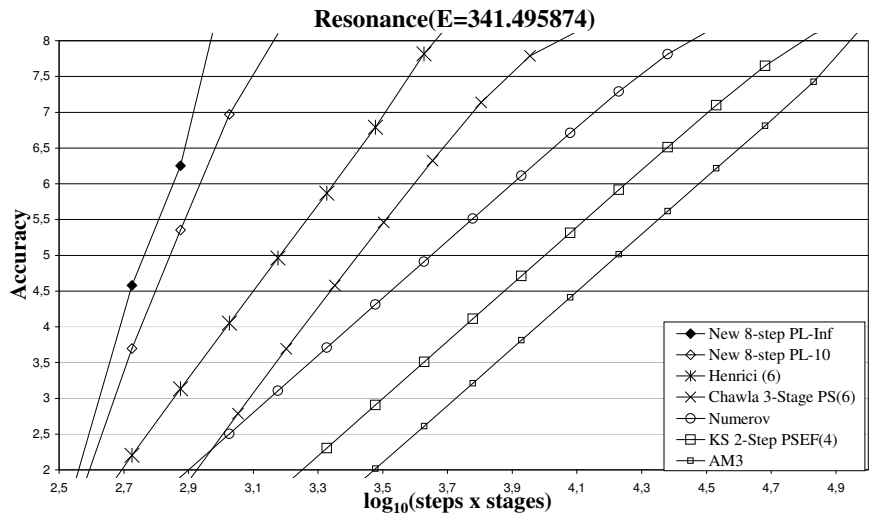}
    \caption{Efficiency for the Schr\"odinger equation using E = 341.495874}
    \label{fig_resonance_341}
\end{figure}

\begin{figure}[htbp]
    \includegraphics[width=\textwidth]{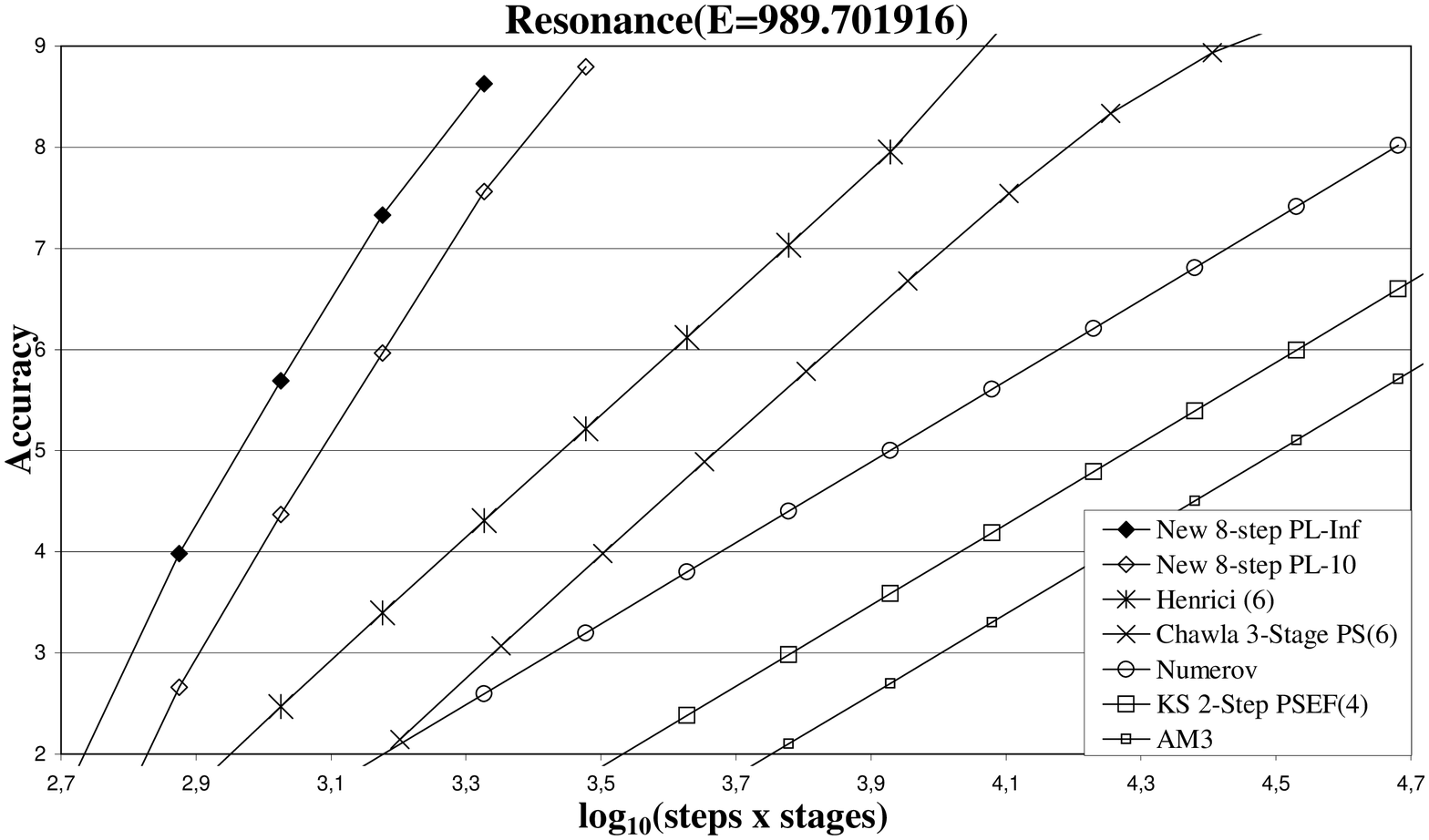}
    \caption{Efficiency for the Schr\"odinger equation using E = 989.701916}
    \label{fig_resonance_989}
\end{figure}

We present the \textbf{accuracy} of the tested methods expressed by the $-\log_{10}$(max. error over interval) or
$-\log_{10}$(error at the end point), depending on whether we know the theoretical solution or not, versus the
$\log_{10}$(steps x stages). In Figure \ref{fig_Orbital Problem by Franco and Palacios} we see the results for the Franco-Palacios almost periodic problem, in Figure \ref{fig_inhom} the results for the Inhomogeneous equation in Figure \ref{fig_twobody} the results for the Two-body problem and in Figure \ref{fig_Duffing} the results for the Duffing equation.
In Figures \ref{fig_resonance_163},\,\ref{fig_resonance_341} and \ref{fig_resonance_989} we see the results for the
Schr\"odinger equation for energies $E= 163.215341$,\  $E = 341.495874$ and $E= 989.701916$ respectively.

Among all the methods used , the new optimized method with infinite order of phase-lag was the most efficient, with the exception
of the Duffing Equation, had almost identical results with the new method with phase-lag order ten.

The difference from the other methods was about 1.2 decimal
digits better for the Schr\"odinger equation for energy $E = 989.701916$ and about 0.7 d.d. for $E= 163.215341$ and $E = 341.495874$. For the other three problems the difference was enormous, where there was an almost vertical increase in the accuracy compared to the
other methods. There were no case where the efficiency dropped below the efficiency of the others.

As regards the other methods, the one of Henrici was the most efficient, with next the method of Chawla, the method
of Numerov and finally the methods of Kalogiratou-Simos and Adams-Moulton.

\section{Conclusions}
We have constructed two optimized eight-step symmetric implicit methods. The first one has phase-lag of order infinite (phase-fitted). The second one has phase-lag of order ten. We have applied the new methods along with a group of recently developed methods from the literature to the Schr\"odinger equation and related problems. We concluded that the new methods are highly efficient compared to other optimized methods which also reveals the importance of phase-lag when solving ordinary differential equations with oscillating solutions.

\end{article}
\end{document}